\documentclass[a4paper,12pt,reqno]{article}
\usepackage{amsmath,amsthm,amssymb,amscd}
\usepackage{amssymb}
\usepackage{amsthm}
\usepackage{amsmath}
\usepackage{latexsym}
\usepackage{epsfig}
\usepackage{amscd}
\usepackage{graphicx,wrapfig}
\usepackage{subcaption}
\usepackage[thinlines,thiklines]{easybmat}
\usepackage[table]{xcolor}
\usepackage{hyperref}
\hypersetup{
    colorlinks=true,
    linkcolor=blue,
    filecolor=cyan,   
    urlcolor=magenta,
    citecolor=red,
    pdftitle={Overleaf Example},
    pdfpagemode=FullScreen,
    }
\def\natu           {\mathbb N}

\def\real		{\mathbb R}

\def\comp		{\mathbb C}

\def\R		{\cal R}

\def\lra		{\longrightarrow}

\title{Periodicity and Circulant Matrices\\ in the Riordan Array of a Polynomial}
 
\author{Nikolai A. Krylov\\ ~ \\
Siena College, Department of Mathematics\\
515 Loudon Road, Loudonville NY 12211, USA\\ ~ \\
nkrylov@siena.edu}

\date{}

\begin{document}

\newtheorem{theorem}{Theorem}
\newtheorem{lemma}[theorem]{Lemma}
\newtheorem{cor}[theorem]{Corollary}
\newtheorem{conj}[theorem]{Conjecture}
\newtheorem{prop}[theorem]{Proposition}
\newtheorem{question}{Problem}
\theoremstyle{definition}
\newtheorem{definition}[theorem]{Definition}
\newtheorem{example}[theorem]{Example}
\numberwithin{equation}{section}

\maketitle
\begin{abstract}
We consider Riordan arrays $\bigl(1/(1-t^{d+1}), ~ tp(t)\bigr)$. These are infinite lower triangular matrices determined by the 
formal power series $1/(1-t^{d+1})$ and a polynomial $p(t)$ of degree $d$.
Columns of such matrix are eventually periodic sequences with a period of $d + 1$, and circulant matrices are used to 
describe the long term behavior of such periodicity when the column's index grows indefinitely. We also  
discuss some combinatorially interesting sequences that appear through the corresponding A - and Z - sequences 
of such Riordan arrays.
\end{abstract}



\section{Introduction}\label{intro}

Let $R$ be an integral domain with 1. 
Take a polynomial of degree $d$, $p(t) = a_0+a_1 t + \cdots + a_d t^d \in R[t]$ 
such that $a_0 \neq 0$ and consider a formal power series (f.p.s.) 
$$
p(t) + p(t)t^{d+1} + \cdots + p(t)t^{n(d+1)} + \cdots = \frac{a_0+a_1 t + 
\cdots + a_d t^d}{1- t^{d+1}} \in R[[t]].
$$
We will think of this f.p.s. as a generating function of the infinite periodic sequence with repeating 
blocks $\{a_0,a_1,\ldots,a_d\}$. Next we assign to $p(t)$ the infinite lower triangular matrix 
(the Riordan array of $p(t)$)

\begin{equation}
\label{array}
RAp = \bigl(1/(1-t^{d+1}), tp(t)\bigr),
\end{equation}
where the $k$-th column ($k \geq 0$) is given by the coefficients of the f.p.s. 
$$
\frac{\bigl(tp(t)\bigr)^k}{1 - t^{d+1}}.
$$
For example, for the constant non-zero polynomial $p(t) = a$ we have the following matrix
\begin{equation}
\label{example1}
RAa = \left( \frac{1}{1-t}, ta\right) = 
\begin{pmatrix}
1 & 0 & 0 & 0 & \cdots \\
1 & a & 0 & 0 & \cdots \\
1 & a & a^2 & 0 & \cdots \\
1 & a & a^2 & a^3 & \cdots \\
 \vdots & \vdots & \vdots & \vdots & \ddots \\
\end{pmatrix}.
\end{equation}

Such matrices are elements of the Riordan group $\R$, which consists of the pairs $(f(t),g(t))$, where 
$f(t) = \sum_{i\geq 0} f_it^i$ and  $g(t) \sum_{i\geq 0} g_it^i $ are certain 
f.p.s. To guarantee nonzero elements along 
the main diagonal and the lower triangular form of the array, the order of $f(t)$ must be 0, that is $f_0 \neq 0$,  
and the order of $g(t)$ must be 1, that is $g_0 = 0 \wedge g_1\neq 0$. The first 
appearance of this group goes back to 1991 when Shapiro, Getu, Woan, and Woodson published their paper 
\cite{Shapiro1}, where they discuss several examples of the arrays with properties similar to those of the Pascal's matrix 
$\bigl(1/(1-t),t/(1-t)\bigr)$ (see also \cite{Sprugnoli1}).

Since then, Riordan matrices have been extensively studied not only as combinatorial objects but also as 
objects with interesting algebraic structure and numerous applications. In this article we will use
only the fundamental properties of such matrices leaving the group structure of $\R$ aside.
For the details, the reader is referred to either of the books 
\cite{Barry} or \cite{Shapiro2}, or a recent survey \cite{Cameron} focusing on combinatorial relationships 
between Riordan matrices and lattice path enumerations. The monograph \cite{Shapiro2} is a 
great resource on recent developments of the topic and the corresponding literature. 

Returning to our matrix $RAa$ in (\ref{example1}), note that we can write the generating function for the 
$k$-th column (starting with $k=0$) as 
$$
a^k t^k + a^k t^{k+1} + \cdots = \frac{a^k t^k}{1-t}.
$$
If we write this $k$-th column as a sequence of numbers 
$$
C_{i,k} = (0,0,\ldots, 0,a^k,a^k,\ldots),
$$
we see that this sequence is {\sl eventually periodic} with period 1 starting with $i = k$, for $i\geq 0$. 
In the next section we will prove that for any polynomial $p(t) \in R[t]$ of degree $d$ each column of the matrix (\ref{array}) 
is eventually periodic with a period of $d+1$. Then in section 3 we will associate to this Riordan array 
$RAp$ a circulant $(d+1) \times (d+1)$ matrix and describe the end behavior 
of the column generating functions of $RAp$ in terms of this matrix and its orbit. We will focus on 
the orbits of linear and quadratic polynomials in section 4, and in the last section, we will 
study the characteristic A- and Z- sequences of $RAp$ and their combinatorial properties.


\section{Column Generating Functions of 
$$
\bigl(1/(1-t^{d+1}), tp(t)\bigr)
$$}

We start with a polynomial $p(t) = a_0 + a_1t + a_2t^2 + \cdots +a_dt^d\in R[t]$ where $a_0 a_d\neq 0$ and 
consider the Riordan array 
$$
RAp = \left( \frac{1}{(1-t^{d+1})}, tp(t)\right).
$$
If $a_0=0$, we will still have eventually periodic columns, but then the main diagonal of our $RAp$ will begin 
with 1 and have 0 after that, so the matrix will not be a {\sl proper Riordan array}. 
If we use sequence $\{e_i\}_{i=0}^{\infty} = (0,\ldots,0,1,0,\ldots ) $ with the only nonzero element when $i=k$,
then by the Fundamental Theorem of Riordan Arrays (see \S 5.1 of \cite{Barry} or 
Theorem 3.1 of \cite{Shapiro2}) the $k$-th column of $RAp$ is 
$$
RAp \cdot (0,\ldots,0,1,0,\ldots )^T  = \left( \frac{1}{(1-t^{d+1})}, tp(t)\right)\cdot t^k = 
\frac{(tp(t))^k}{(1-t^{d+1})}.
$$

\begin{theorem} For every integer $k\geq 0$,
$$
\frac{(tp(t))^k}{(1-t^{d+1})}
$$
is the generating function of an eventually periodic sequence $\{C_{i,k}\}_{i=0}^{\infty}$ with a period of $d + 1$.
The periodicity starts at the $1+ (k-1)(d+1)$-st place, i.e. for any $n\geq 0$ we have
\begin{equation}
\label{formula1}
C_{1+(k-1)(d+1) + n, k} = C_{1+ k(d+1) + n, k}.
\end{equation}
\end{theorem}
\begin{proof}
We use induction on $k$. Since the zeroth column of $RAp$ is given by the geometric series 
$1/(1-t^{d+1})$, the base of induction is obvious. Now let's assume that the statement is true for a $k$-th column 
($k\geq 0$). Then $k+1$-st column of $RAp$ is given by a f.p.s. 
$$
\sum\limits_{i=0}^{\infty} C_{i,k+1}t^i= \frac{(tp(t))^{k+1}}{(1-t^{d+1})} = \frac{(tp(t))^k}{(1-t^{d+1})}\cdot tp(t),
$$
and using the notation $[t^n]f(t)$ for the {\sl extraction of the coefficient of $t^n$ from a f.p.s. $f(t)$}, 
we must show that for every integer $n\geq 0$,
$$
[t^{1+ k(d+1) + n}]\frac{(tp(t))^{k+1}}{(1-t^{d+1})} = 
[t^{1+(k+1)(d+1) + n}]\frac{(tp(t))^{k+1}}{(1-t^{d+1})}.
$$
Since $\deg(tp(t)) = d+1$ the convolution rule for the coefficient operator (see for example section 
2.2 of \cite{Shapiro2}) for all $n\geq 0$ gives
$$
C_{1+(k+1)(d+1) + n, k+1} = [t^{1+(k+1)(d+1) + n}]\frac{(tp(t))^k}{(1-t^{d+1})}\cdot (tp(t)) 
$$

$$
= \sum\limits_{j=0}^{1+(k+1)(d+1) + n}[t^j]\bigl(tp(t)\bigr) \cdot [t^{1+(k+1)(d+1) + n - j}]\frac{(tp(t))^k}{(1-t^{d+1})} 
$$

\begin{equation}
\label{Thm1_1}
= \sum\limits_{j=1} ^ {d+1} a_{j-1} [t^{1+ k(d+1) + (n + d + 1 - j)}]\frac{(tp(t))^k}{(1-t^{d+1})}.
\end{equation}
Using induction hypothesis for all $n\geq 0$ and $0\leq j\leq d+1$ we have 
$$
C_{1+ k(d+1) + (n + d + 1 - j), k} = C_{1+ (k-1)(d+1) + (n + d + 1 - j), k},
$$
that is  
$$
[t^{1+ k(d+1) + (n + d + 1 - j)}]\frac{(tp(t))^k}{(1-t^{d+1})} = [t^{1+ (k-1)(d+1) + (n +d+1 - j)}]\frac{(tp(t))^k}{(1-t^{d+1})}
$$
and we can continue (\ref{Thm1_1}) as 
$$
\sum\limits_{j=1} ^ {d+1} a_{j-1} \cdot [t^{1+ (k-1)(d+1) +( n +d+1 - j)}]\frac{(tp(t))^k}{(1-t^{d+1})} = 
\sum\limits_{j=1} ^ {d+1} a_{j-1} \cdot [t^{1+ k(d+1) + n - j}]\frac{(tp(t))^k}{(1-t^{d+1})}
$$

$$
= \sum\limits_{j=0}^{1+ k(d+1) + n}[t^j]\bigl(tp(t) \bigr) \cdot [t^{1+ k(d+1) + n - j}]\frac{(tp(t))^k}{(1-t^{d+1})} = 
C_{1+ k(d+1) + n, k+1},
$$
which finishes induction and proves (\ref{formula1}).
\end{proof}

Here is an example to illustrate this theorem with $p(t) = 1 + 5t$:
\begin{equation}
\label{matrix1}
\left( \frac{1}{1- t^2}, t + 5t^2\right) = 
\left(
\begin{BMAT}[3pt]{ccccccc}{cccccccc}
1 & 0 & 0 & 0 & 0 & 0 & 0  \\
0 & 1 & 0 & 0 & 0 & 0 & 0  \\
1 & 5 & 1 & 0 & 0 & 0 & 0  \\
0 & 1 &10 & 1 & 0 & 0 & 0  \\
1 & 5 & 26 & 15 & 1 & 0 & 0  \\
0 & 1 & 10 & 76 & 20 & 1 & 0  \\
1 & 5 & 26 & 140 & 151 & 25 & 1 \\
\vdots & \vdots &  \vdots & \vdots & \vdots & \vdots & \ddots 
\addpath{(1,5,1)ruuldd}
\addpath{(2,3,1)ruuldd}
\addpath{(3,1,1)ruuldd}
\end{BMAT}
\right)
\end{equation}
The third column begins with $(0,0,0,1,15,76,140,76,140,\ldots)$, and periodicity 
starts with 76, whose index is $i = 5$. Also note that 

\begin{equation}
\label{orbit1}
\begin{pmatrix}
5 & 1\\
1 & 5
\end{pmatrix}
\begin{bmatrix}
1\\
5
\end{bmatrix} = 
\begin{bmatrix}
10\\
26
\end{bmatrix} ~ \mbox{and} ~ 
\begin{pmatrix}
5 & 1\\
1 & 5
\end{pmatrix}
\begin{bmatrix}
10\\
26
\end{bmatrix} = 
\begin{pmatrix}
5 & 1\\
1 & 5
\end{pmatrix}^2
\begin{bmatrix}
1\\
5
\end{bmatrix} = 
\begin{bmatrix}
76\\
140
\end{bmatrix},
\end{equation}
where the columns of the $2\times 2$ matrix are obtained from the coefficients of 
$p(t) = 1 + 5t$ by a cyclic permutation. Such square matrices belong to a special type of Toeplitz matrices, 
called {\sl circulant} matrices. We will see in the next section that the periodic 
cycles in the columns of our Riordan array make the orbit of the 
coefficient sequence of $p(t)$ under the linear map defined by such a matrix.


\section{Circulant Matrix of $p(t)$ and its Orbit}

For a polynomial $p(t) = a_0+a_1 t + \cdots + a_d t^d\in R[t]$ of degree $d$ with $a_0 \neq 0$ let us 
write its coefficients as a row vector starting with $a_d$ 
$$
\nu = (a_d, a_{d-1}, \ldots ,a_0)\in R^{d+1},
$$
and consider a shift operator $T:R^{d+1}\to R^{d+1}$ defined as 
$$
T(a_d, a_{d-1},\ldots, a_1,a_0):= (a_0, a_d, a_{d-1},\ldots, a_1).
$$
Following \cite{Kra}, we associate to the vector $\nu$ the $(d+1)\times (d+1)$ matrix $V_{p(t)} = {\rm circ} \{\nu\}$, 
and call it the {\sl circulant matrix} associated to the polynomial $p(t)$. 
Its rows are given by iterations of the shift operator acting on the vector $\nu$, i.e. the
$i$-th row of $V_{p(t)}$ is $T^i \nu,~ i\in\{0,\ldots,d \}$:
$$
V_{p(t)} = 
\begin{pmatrix}
a_d & a_{d-1} & \cdots & a_1 & a_0\\
a_0 & a_{d} & \cdots & a_2 & a_1\\
\vdots & \vdots & \ddots & \vdots & \vdots\\
a_{d-2} & a_{d-3} & \cdots & a_d & a_{d-1}\\
a_{d-1} & a_{d-2} & \cdots & a_0 & a_d\\
\end{pmatrix}.
$$
We refer the reader to the monographs \cite{Davis} and \cite{Fuhrmann} for the detailed 
treatment of the subject. Denote the $n$-th primitive root of unity by $\xi = e^{2\pi i/n}$ and for each 
$l\in \{0,\ldots, n-1\}$ consider vector $x_l = \frac{1}{\sqrt{n}}(1,\xi^l,\xi ^{2l},\ldots, \xi^{(n-1)l})\in \comp^n$. 
Using $x_l$ as rows define the Fourier matrix
$$
F:=\begin{pmatrix} x_0\\ x_1 \\ \vdots \\ x_{n-2} \\ x_{n-1} \end{pmatrix} = \frac{1}{\sqrt{n}} \begin{pmatrix}
1 & 1 & \ldots & 1 & 1\\
1 & \xi  & \ldots & \xi^{n-2} & \xi^{n-1} \\
\vdots & \vdots  & \ddots & \vdots & \vdots \\
1 & \xi^{n-2}  & \ldots & \xi^{(n-2)(n-2)} & \xi^{(n-2)(n-1)}\\
1 & \xi^{n-1} &  \ldots & \xi^{(n-1)(n-2)} & \xi^{(n-1)(n-1)}\\
\end{pmatrix}.
$$
This is the matrix of the discrete (or finite) Fourier transform (\cite{Davis}, \S 2.5), which is unitary, has order four, and satisfies 
$$
F = F^T, ~~ F^* = (F^*)^T = \overline{F}, ~~  F = \overline{F}^*, ~~ F\cdot F^* = F^*\cdot F = I.
$$
It turns out that for any circulant matrix $V$ of size $n\times n$, $V$ is diagonalized by $F$, that is
\begin{equation}
\label{diag}
F^{-1} V F = D_V,
\end{equation}
where $D_V$ is the diagonal matrix (Theorem 12. of \cite{Kra}). Moreover, the eigenvalues of 
$V: \lambda_0,\lambda_1,\ldots,\lambda_{n-1}$ are the values of the {\sl representer polynomial} $P_V$ 
at $\xi^l,~l\in\{0,1,\ldots, n-1\}$ (Theorem 6 of \cite{Kra} or Theorem 3.2.2 of \cite{Davis}). The representer 
polynomial of $V_{p(t)}$ is defined by taking the dot product of the vector $\nu$ 
with $(1,Z,\ldots,Z^d)$, that is 
$$
P_{V_{p(t)}}(Z) = \sum\limits_{i=0}^{d} a_{d-i} Z^i = Z^d\cdot p(1/Z).
$$
Thus, as mentioned above, we have
$$
D_V = {\rm diag}\bigl(P_V(1),P_V(\xi),\ldots,P_V(\xi^{n-1})\bigr).
$$
Now let us confirm the relationship, which will generalize equalities (\ref{orbit1}). Using the same notations 
as in the Theorem 1 of section 2 for the generating functions of the columns of $RAp$, we have the following

\begin{theorem}
Let $p(t) = a_0+a_1 t + \cdots + a_d t^d = C_{1,1} + C_{2,1} t + \cdots + C_{1+d,1} t^d$ be as above. 
Then for all $n\geq 0$, the periodic block of $(n+1)$-st column of (\ref{array}) is the $n$-th iteration of $V_{p(t)}$ 
applied to the periodic block of the first column, that is
\begin{equation}
\label{orbit}
V_{p(t)}^n \begin{pmatrix} a_0 \\ a_1 \\ \vdots \\ a_d \end{pmatrix} =
\begin{pmatrix}
C_{1+n(d+1), n+1} \\ C_{1+n(d+1)+1, n+1}  \\ \vdots \\ C_{1+n(d+1)+d, n+1} \\
\end{pmatrix}.
\end{equation}
\end{theorem}
\begin{proof}
The base of induction with $n = 0$ is obvious, so we assume that (\ref{orbit}) is true for any $n$ 
and show that (\ref{orbit}) will remain true when $n$ is replaced by $n+1$. Let us look at the 
$(1+(n+1)(d+1)+k)$-th coefficient of the $(n+2)$-nd column represented by the f.p.s. as
$$
C_{1+(n+1)(d+1)+k, n+2} = [t^{1+(n+1)(d+1)+k}] \left(tp(t) \frac{(tp(t))^{n+1}}{1-t^{d+1}} \right).
$$

If we, for a moment, abbreviate $1+n(d+1)$ by $H$, and use the convolution formula for the product of two 
series, we obtain
$$
C_{1+(n+1)(d+1)+k, n+2} = \sum\limits_{j = 0}^{1+(n+1)(d+1)+k} [t^j](tp(t)) 
[t^{1+(n+1)(d+1)+k-j}]\left(\frac{(tp(t))^{n+1}}{1-t^{d+1}} \right)
$$
\begin{equation}
\label{sum0}
= \sum\limits_{j=1}^{d+1} a_{j-1}C_{1+n(d+1)+k+(d+1)-j, n+1} = 
\sum\limits_{j=1}^{d+1} a_{j-1}C_{H+k+(d+1)-j, n+1}.
\end{equation}
To simplify the notations, let me also drop the column index $n+1$. Then (\ref{sum0}) equals 
\begin{equation}
\label{sum12}
a_0C_{H+(d+1)+k-1} + a_1C_{H+(d+1)+k-2} + \ldots + a_{k-2}C_{H+(d+1)+1} + 
a_{k-1}C_{H+(d+1)}  
\end{equation}
\begin{equation}
\label{sum2}
+ a_{k}C_{H+(d+1)-1} + \cdots + a_dC_{H+k}.
\end{equation}
Using periodicity (\ref{formula1}) for every $s\in\{k,\ldots, d\}$ we have
$
C_{H+s} = C_{H+(d+1)+s},
$
so (\ref{sum2}) can be written as
\begin{equation}
\label{sum14}
\sum\limits_{i = k+1}^{d+1} a_{d+1 + k - i}C_{1+(n+1)(d+1) + i - 1}.
\end{equation}
Therefore, if we read (\ref{sum12}) from right to left and then add to it (\ref{sum14}), the result 
will be equal to the summation 
$$ 
\sum\limits_{i =1}^{d+1} \widehat{a_{k-i}}C_{1+(n+1)(d+1) + i - 1},
$$
where $\widehat{a_{k-i}} \in \{0,\ldots,d\}$  and $\widehat{a_{k-i}} \equiv a_{k-i} \pmod{d+1}$.
But (\ref{sum12}) together with (\ref{sum14}) equals (\ref{sum0}), hence 
$$
C_{1+(n+1)(d+1)+k, n+2} = \sum\limits_{i =1}^{d+1} \widehat{a_{k-i}}C_{1+(n+1)(d+1) + i - 1, n + 1}.
$$
This shows that for every $k\in\{0,\ldots, d\}$, the 
$(1+(n+1)(d+1)+k)$-th coefficient of the $(n+2)$-nd column of the Riordan array $RAp$ equals the 
product of the corresponding row of $V_{p(t)}$ and the periodic block of the $(n+1)$-st column. That is 
$$
C_{1+(n+1)(d+1)+k, n+2} = \bigl(T^{k}(a_d, a_{d-1},\ldots, a_1,a_0)\bigr)\cdot 
\begin{pmatrix}
C_{1+n(d+1), n+1} \\ C_{1+n(d+1)+1, n+1}  \\ \vdots \\ C_{1+n(d+1)+d, n+1} \\
\end{pmatrix},
$$
which proves the induction step. 
\end{proof}

\begin{cor}
Columns of the Riordan array $RAp$ of $p(t) = a_0 + a_1t + \cdots + a_dt^d$ contain the complete 
forward orbit of the coefficient vector $(a_0,a_1,\ldots,a_d)^T$ under the iterations of the 
circulant matrix $V_{p(t)}$. 
\end{cor}


\section{Orbits of Linear and\\ Quadratic Polynomials.}

In the rest of the paper we will consider polynomials and the formal power series over $\real$. 
First we give an alternative description of the orbit given by (\ref{orbit}) using equality (\ref{diag}).
Notice that for $p(t) = a_0 + a_1t + \cdots + a_dt^d$, and $\xi = e^{2\pi i/(d+1)}$ with 
the corresponding Fourier matrix $F$, we have the following identities (cf. (2.5.10) and (2.5.13) of \cite{Davis}):
$$
\sqrt{d+1} F 
\begin{pmatrix}
a_0\\a_1\\ \vdots \\ a_d\\ 
\end{pmatrix} = 
\begin{pmatrix}
p(1) \\ p(\xi) \\ \vdots \\ p(\xi^d) \\
\end{pmatrix} ~~~~ \mbox{and} ~~~~ 
\sqrt{d+1} \cdot \overline{F} 
\begin{pmatrix}
a_0\\a_1\\ \vdots \\ a_d\\ 
\end{pmatrix} = 
\begin{pmatrix}
p(1) \\ p(\overline{\xi}) \\ \vdots \\ p(\overline{\xi^d}) \\
\end{pmatrix}.
$$
Since 
$$
V_{p(t)}^n = F D_{V_{p(t)}}^n \overline{F},
$$
and the diagonal matrix here is $D_V = {\rm diag}\bigl(P_V(1),P_V(\xi),\ldots,P_V(\xi^{d})\bigr)$, with 
$P_V(\xi^j) = \xi^{jd}\cdot p(\xi^{-j})$, we deduce that for all $n\in\natu$ 

$$
V_{p(t)}^n \begin{pmatrix} a_0 \\ \vdots \\ a_j \\ \vdots \\ a_d \end{pmatrix} = \frac{1}{\sqrt{d+1}}  F 
\begin{pmatrix}
p^n(1) \cdot p(1) \\ \vdots \\ \xi^{njd} p^n(\xi^{-j})\cdot p(\overline{\xi}^j) \\ \vdots \\ 
\xi^{nd^2}p^n(\xi^{-d})\cdot p(\overline{\xi}^d) \\
\end{pmatrix} = \frac{1}{\sqrt{d+1}} F 
\begin{pmatrix}
p^{n+1}(1) \\ \vdots \\ \xi^{njd} p^{n+1}(\xi^{-j})\\ \vdots\\ \xi^{nd^2}p^{n+1}(\xi^{-d}) \\
\end{pmatrix}.
$$
The last equality follows from the fact that for any root of unity $\xi$, the conjugate $\overline{\xi}$ and the 
multiplicative inverse $1/{\xi}$ are equal. Multiplying out the matrix $F$ by the column vector 
$\left( \ldots, \xi^{njd}p^{n+1}(\xi^{-j}),\ldots\right)^T$ we obtain an alternative presentation of the orbit as

\begin{equation}
\label{orbit2}
\frac{1}{d+1}\begin{pmatrix}
\sum\limits_{k=0}^{d}\xi^{nkd} p^{n+1}(\xi^{-k})  \\

\vdots \\

\sum\limits_{k=0}^{d}\xi^{kj + nkd} p^{n+1}(\xi^{-k})  \\
 
\vdots\\

\sum\limits_{k=0}^{d}\xi^{(n+1)kd} p^{n+1}(\xi^{-k}) \\
\end{pmatrix}.
\end{equation}

Recall that for every $k\geq 0$, the $k$th eigenvalue of $V_{p(t)}$ is $\lambda_{k} = \xi^{dk}p(\xi^{-k})$, 
hence it is clear that we can rewrite (\ref{orbit2}) as a product of a matrix by the eigenvalues vector 
of $V_{p(t)}$. We state this description as 
\begin{theorem}
The forward orbit of $(a_0,a_1,\ldots,a_d)^T$ under $V_{p(t)}$ is the sequence 
$$
V_{p(t)}^n \begin{pmatrix} a_0 \\ a_1 \\ \vdots \\ a_{d-1} \\ a_d \end{pmatrix} = \frac{1}{d+1} 
\begin{pmatrix}
1 & \xi & \cdots & \xi^{d-1} & \xi^d\\
1 & \xi^2 & \cdots & \xi^{2(d-1)} & \xi^{2d} \\
\vdots & \vdots & \ddots & \vdots & \vdots \\
1 & \xi^d & \cdots & \xi^{(d-1)d} & \xi^{d^2}\\
1 & 1 & \cdots & 1 & 1\\
\end{pmatrix}
.\begin{pmatrix}
\lambda^{n+1}_0 \\
\lambda^{n+1}_1 \\
\vdots \\
\lambda^{n+1}_{d-1} \\
\lambda^{n+1}_{d} \\
\end{pmatrix}
$$

\end{theorem}
\begin{proof}
It follows directly from (\ref{orbit2}) using $\xi^{d+1} = 1$.
\end{proof}

Notice that the matrix in this theorem can be obtained from the Fourier matrix $F$ by a cyclic 
permutation that shifts the first row to the bottom and pushes the others one level up.

Before discussing specifically orbits of linear and quadratic polynomials, allow me alter our point of view here. 
It is obvious from the definition of $V_{p(t)}$ that $(1,1,\ldots,1)^T\in\real^{d+1}$ is the eigenvector 
corresponding to the eigenvalue $\lambda_0 = \xi^0 p(\xi^{-0}) = \sum\limits_{i=0}^{d} a_i$. To simplify the 
notations and analysis, let us rotate $\real^{d+1}$ so that $(0,\ldots,0,\sqrt{d+1})^T\in\real^{d+1}$
maps to $(1,1,\ldots,1)^T\in\real^{d+1}$. We choose for such a rotation 
an orthogonal matrix $R$ of determinant 1. Such transformation preserves the inner product in 
$\real^{d+1}$, and asymptotic behavior of the orbit clearly doesn't depend on the choice of a rotation. 
In the new basis the orbit of the coefficient sequence of $p(t)$ will be given by the sequence
\begin{equation}
\label{rotation}
\bigl(R \cdot \widetilde{V}_{p(t)}\cdot R^{-1} \bigr)^n
\begin{pmatrix} a_0 \\ a_1 \\ \vdots \\ a_d \end{pmatrix}
= R \cdot \widetilde{V}^n_{p(t)}\cdot R^{-1} 
\begin{pmatrix} a_0 \\ a_1 \\ \vdots \\ a_d \end{pmatrix},
\end{equation}
where $\widetilde{V}_{a+bt} = R^{-1} V_{a+bt} R$, $\forall n\in\natu$.
It will be easier to describe the long term behavior of our orbits if we apply such rotations in $\real^2$ and 
$\real^3$.


\subsection{Orbit of a Linear Polynomial $p(t) = a + bt,~a\cdot b\neq 0$}

Since $\xi = e^{2\pi i/2} = -1$ and $p^{n+1}(\pm 1) = (a \pm b)^{n+1}$, we have $\lambda_0 = (a+b)$ and 
$\lambda_1 = (-1)(a-b)$. Hence the orbit of $p(t) = a + bt$ is given by 
\begin{equation}
\label{rotation1}
\begin{pmatrix}
b & a\\
a & b\\
\end{pmatrix}^n 
\begin{pmatrix} 
a \\ b \\
\end{pmatrix}  = \frac{1}{2}
\begin{pmatrix}
1 & -1\\
1 & 1
\end{pmatrix}
\begin{pmatrix}
(a+b)^{n+1}\\
(b-a)^{n+1}
\end{pmatrix} = \frac{1}{2}
\begin{pmatrix}
(a+b)^{n+1} - (b-a)^{n+1} \\
(a+b)^{n+1} + (b-a)^{n+1} \\
\end{pmatrix}.
\end{equation}
For the following orthogonal matrix $R$, which rotates $\real^2$ clockwise by $\pi/4$, 
$$
R = \begin{pmatrix}
1/\sqrt{2} & 1/\sqrt{2}\\
-1/\sqrt{2} & 1/\sqrt{2}
\end{pmatrix},
$$
we can use $R^{-1} V^n_{a+bt} = \widetilde{V}^n_{a+bt} R^{-1}$ and write the {\sl rotated orbit} as 
$$
R^{-1} V_{a+bt} ^n 
\begin{pmatrix} 
a \\ b \\
\end{pmatrix}
= \widetilde{V}^n_{a+bt}
\begin{pmatrix} 
(a - b)/\sqrt{2} \\ (a + b)/\sqrt{2} \\
\end{pmatrix}
$$
\begin{equation}
\label{rotation2} 
= \begin{pmatrix}
b - a & 0\\
0 &  a + b
\end{pmatrix}^n 
\begin{pmatrix} 
(a - b)/\sqrt{2} \\ (a + b)/\sqrt{2} \\
\end{pmatrix} = 
\begin{pmatrix}
-(b - a)^{n+1}/\sqrt{2} \\ 
(a + b)^{n+1}/\sqrt{2} \\ 
\end{pmatrix}.
\end{equation}
Formulas (\ref{rotation2}) and (\ref{rotation1}) are of course equivalent, but the second one is easier 
to work with. 

The origin $(0,0)$ is the fixed point of $V_{a+bt}$, as well as of $\widetilde{V}_{a+bt}$, unless one of the 
eigenvalues is 1. So let us discuss first the case when $a + b = 1$, the case when $b - a = 1$ is treated similarly. 
Assuming that $|b - a| = |1 - 2a| < 1$, we see that the orbit in the standard basis converges to the point 
$(1/2,1/2)$, and in the rotated one to $(0,1/\sqrt{2})$. When we rotate $(0,1/\sqrt{2})$ clockwise by 
$\pi/4$ we will get the point $(1/2,1/2)$ as expected, so we will focus from now on 
the orbit given by (\ref{rotation2}). When $a + b = 1$ and $|b-a| = |1 - 2a| > 1$, the orbit converges to
infinity staying on the line $\tilde{y} = 1/\sqrt{2}$ (in the rotated coordinates). If $b - a < -1$, it splits into 
two parts (one diverging to positive infinity and the other one to the negative infinity) accordingly 
if $n$ is odd or even. Now, if both eigenvalues are less than 1 in absolute value, the orbit contracts to the 
fixed point $(0,0)$ along the curve given parametrically (again, in the rotated coordinates) as 
$$
\tilde{x}(t) = \frac{-1}{\sqrt{2}}(b-a)^t, ~ \tilde{y}(t) = \frac{1}{\sqrt{2}}(a+b)^t,
$$
or (compare with \S 3.1 of \cite{Katok}) by 
$$
|\tilde{y}| = C|\tilde{x}|^{\ln|a + b|/\ln|b - a|} ~~ \mbox{for some constant} ~ C.
$$
For example, if $a = -4/11$ and $b = 6/11$, then 
$$
C = \sqrt{2}^{\left(\frac{\ln (2/11)}{\ln (10/11)} -1\right)} 
\approx 348.05187 ~~ \mbox{and} ~~ y = C (-x)^{\frac{\ln(2/11)}{\ln(10/11)}}.
$$ 
Figure 1 here shows the curve together with a first few iterations of the point $(-4/11,6/11)$.

\begin{figure}[h] 
\centering
\includegraphics[width=70mm]{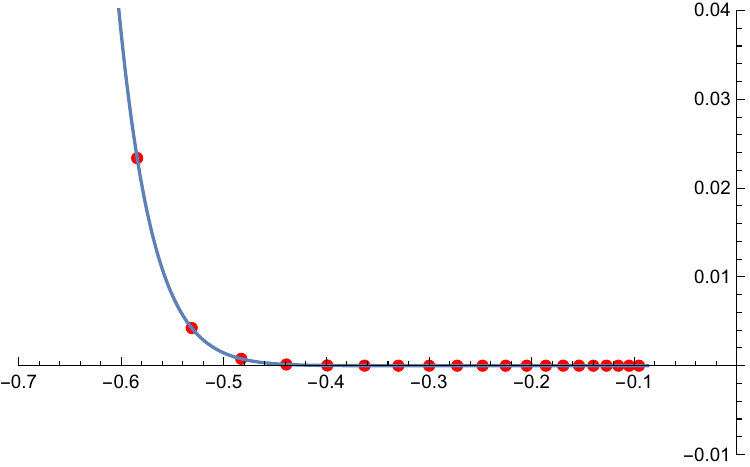}
\caption{Orbit (\ref{rotation2}) of  $(-4/11,6/11)$}
\label{}
\end{figure}

On the other hand, when $|a+b|>1$ and $|a-b|>1$, the orbit clearly tends to infinity exponentially. 
I will end the case of $p(t) = a + bt$ by noting that the division of the orbit into two different 
subsequences can also happen when one of the 
eigenvalues equals $-1$ and the absolute value of the other one is less than 1. For example, if 
$\lambda_0 = a+b= -1/m, m > 1$ and $\lambda_1 = b-a = -1$, the orbit splits into two convergent 
subsequences 
$$
\begin{pmatrix}
-(-1)^{2k} /\sqrt{2} \\
(-1/m)^{2k}/ \sqrt{2}\\
\end{pmatrix} \stackrel{k\to\infty}{\lra}
\begin{pmatrix}
-1/\sqrt{2} \\ 0 \\
\end{pmatrix} ~~ \mbox{and} ~~ 
\begin{pmatrix}
-(-1)^{2k+1} /\sqrt{2} \\
(-1/m)^{2k+1} /\sqrt{2}\\
\end{pmatrix} \stackrel{k\to\infty}{\lra}
\begin{pmatrix}
1/\sqrt{2} \\ 0 \\
\end{pmatrix}.
$$


\subsection{Orbit of a Quadratic $p(t) = a + bt + ct^2$}

Let us focus on the orbit given by $R^{-1} V^n_{p(t)} = \widetilde{V}^n_{p(t)} R^{-1}$ and 
choose the rotation matrix to be 
$$
R = \begin{pmatrix}
1/\sqrt{6} & -1/\sqrt{2} & 1/\sqrt{3} \\
1/\sqrt{6} & 1/\sqrt{2} & 1/\sqrt{3} \\
-\sqrt{2/3} & 0 &  1/\sqrt{3} \\
\end{pmatrix}.
$$
Since $\xi = e^{2\pi i/3}  = \cos(2\pi/3) + i \sin(2\pi/3) = (-1+\sqrt{-3})/2$ and $Z^2p(1/Z) = aZ^2 + bZ + c$, 
we have $\lambda_0 = a + b + c $, $\lambda_1 = a\xi^2+b\xi + c $, and 
$\lambda_2 = a\xi+b\xi^2 + c $. The orbit will be given by the sequence 
\begin{equation}
\label{rotation3}
\widetilde{V}^n_{p(t)} R^{-1}
\begin{pmatrix}
a \\ b\\ c \\
\end{pmatrix} =
\begin{pmatrix}
\frac{2c - a - b}{2} & \frac{(b - a)\sqrt{3}}{2} & 0 \\
\frac{(a - b)\sqrt{3}}{2} & \frac{2c - a - b}{2} & 0 \\
 0 & 0 & a + b + c\\
\end{pmatrix}^n 
\begin{pmatrix}
\frac{a+b-2c}{\sqrt{6}} \\ \frac{b-a}{\sqrt{2}} \\  \frac{a + b + c}{\sqrt{3}} \\
\end{pmatrix}.
\end{equation}

The matrix $\widetilde{V}_{p(t)} =  R^{-1} V_{p(t)} R$ in (\ref{rotation3}) is a scaling-rotation matrix in $\real^3$, 
and since $\det(R) = 1$, 
$$
\det(\widetilde{V}_{p(t)}) = \det(V_{p(t)})  = (a+b+c)(a^2+b^2+c^2 - ab - ac - bc).
$$
The scaling factor along 
the third basis vector is obviously $(a+b+c)$. Since $\forall a,b\in\real$, $a^2 + b^2 \geq 2ab$, it follows that 
$a^2+b^2+c^2-(ab+ac+bc) \geq 0$ and $|2c - a - b| \leq 2\sqrt{a^2+b^2+c^2-ab-ac-bc}$, $\forall a,b,c\in\real$.
Therefore, if we exclude the exceptional (but trivial) case when $a = b = c$, we can rewrite 
the top $2\times 2$ block matrix as 
$$
\sqrt{a^2+b^2+c^2-(ab+ac+bc)}
\begin{pmatrix}
\cos\theta & \sin\theta \\
-\sin \theta  & \cos \theta \\
\end{pmatrix},
$$
where
\begin{equation}
\label{cos}
\cos \theta = \frac{2c - a - b }{2\sqrt{a^2+b^2+c^2-ab-ac-bc}}.
\end{equation}
When $2c < a+b $, the fraction in (\ref{cos}) is negative and the orbit splits naturally into two 
opposite (i.e. rotated by $\pi$) spirals. It happens because when values of $n$ are odd, we have
$$
\begin{pmatrix}
\cos\theta  & \sin\theta \\
-\sin \theta  & \cos \theta \\
\end{pmatrix}^n = 
\begin{pmatrix}
-|\cos\theta | & \sin\theta \\
-\sin \theta  & - |\cos \theta | \\
\end{pmatrix}^n = (-1) 
\begin{pmatrix}
|\cos\theta | & - \sin\theta \\
\sin \theta  & |\cos \theta | \\
\end{pmatrix}^n,
$$
but for the even values of $n$,

$$
\begin{pmatrix}
\cos\theta  & \sin\theta \\
-\sin \theta  & \cos \theta \\
\end{pmatrix}^n = 
\begin{pmatrix}
-|\cos\theta | & \sin\theta \\
-\sin \theta  & - |\cos \theta | \\
\end{pmatrix}^n = 
\begin{pmatrix}
|\cos\theta | & - \sin\theta \\
\sin \theta  & |\cos \theta | \\
\end{pmatrix}^n.
$$
Such behavior is similar to the one we saw above  
for a linear polynomial $p(t) = a + bt$ with 
$a + b = 1$ and $b - a < -1$. Thus, depending on whether 
$\sqrt{a^2+b^2+c^2-(ab+ac+bc)}$ is less or greater or equal 1, 
the orbit will either move along a spiral(s) towards 0, expand to infinity, or will stay on a 
circular cylinder $\tilde{x}^2 + \tilde{y}^2 = 1$ correspondingly. I use 
$(\tilde{x},\tilde{y},\tilde{z})$ for the coordinates in the rotated by $R$ basis. 
Below are the parametric equations of the ``helical" 
curve(s) containing the orbit. To simplify the notations I use $r = \sqrt{a^2+b^2+c^2-(ab+ac+bc)}$.
\begin{equation}
\label{helix}
\begin{Bmatrix}
\tilde{x}(t) \\
\tilde{y}(t) \\
\tilde{z}(t)
\end{Bmatrix} = 
\begin{Bmatrix}
\frac{r^t\bigl((a + b - 2c)\cos \theta t + \sqrt{3}(b-a)\sin \theta t\bigr)}{\sqrt{6}}\\
\frac{r^t\bigl((2c - a - b)\sin \theta t + \sqrt{3}(b-a)\cos \theta t\bigr)}{\sqrt{6}} \\
\frac{(a+b+c)^{t+1}}{\sqrt{3}}
\end{Bmatrix},~~ \mbox{when}~\cos \theta \geq 0
\end{equation}
and similarly, if $\cos \theta < 0$ we have the equations for two curves 

$$
\begin{Bmatrix}
\frac{r^t\bigl((a + b - 2c)\cos \theta t + \sqrt{3}(b-a)\sin \theta t\bigr)}{\sqrt{6}}\\
\frac{r^t\bigl((2c - a - b)\sin \theta t + \sqrt{3}(b-a)\cos \theta t\bigr)}{\sqrt{6}} \\
\frac{(a+b+c)^{t+1}}{\sqrt{3}}
\end{Bmatrix}
\bigcup 
\begin{Bmatrix}
\frac{-r^t\bigl((a + b - 2c)\cos \theta t + \sqrt{3}(b-a)\sin \theta t\bigr)}{\sqrt{6}}\\
\frac{-r^t\bigl((2c - a - b)\sin \theta t + \sqrt{3}(b-a)\cos \theta t\bigr)}{\sqrt{6}} \\
\frac{(a+b+c)^{t+1}}{\sqrt{3}}
\end{Bmatrix}.
$$

Figures 2, 3 and 4 on the last page show a few examples of several first iterations of such orbits. 
The corresponding numerical data are given in the table here.

~

\begin{tabular}{||c||c||c||c||}
\hline
\multicolumn{4}{||c||}{Numerical parameters for three different orbits.}\\  
\hline\hline
1. & $(0.93, 0.5, -0.38)$ & $(-0.5, 0.4, 0.89)$ & $(0.9289, 0.487, -0.2159)$\\
\hline\
2.  & 1.05 - escapes & 0.79 - converges & 1.2 - escapes \\
\hline\
3. & $\cos \theta \approx -0.947$ & $\cos \theta \approx 0.77$ & $\cos \theta \approx -0.924$\\
\hline\
4. & $ r\approx 1.15659$ & $ r\approx 1.2211$ & $ r\approx 1$\\
\hline
\end{tabular}

\begin{enumerate}
\item Standard coordinates of the coefficients $(a,b,c)$ of $p(t) = a + bt + ct^2$.

\item Scaling factor $|a + b + c|$: When $|a + b + c| > 1$, the orbit tends to infinity along the third axes $\tilde{z}$, 
when $|a + b + c| < 1$, it converges to the plane $(\tilde{x},\tilde{y},\tilde{0})$.

\item $\cos(\theta)$: When $\cos(\theta) \geq 0$, the orbit is a single spiral (or a circle), 
when $\cos(\theta) < 0$ the orbit consists of two spirals.

\item Scaling factor $r = \sqrt{a^2+b^2+c^2-(ab+ac+bc)}$: When $r < 1$, the orbit contracts to the center of 
the plane $(\tilde{x},\tilde{y},\tilde{0})$, when $r = 1$, it stays on a circular cylinder, and if $r > 1$, the parallel 
projection of the orbit onto the plane $(\tilde{x}, \tilde{y},0)$ grows unbounded.
\end{enumerate}

\begin{figure}[h] 
\centering
\includegraphics[width=65mm]{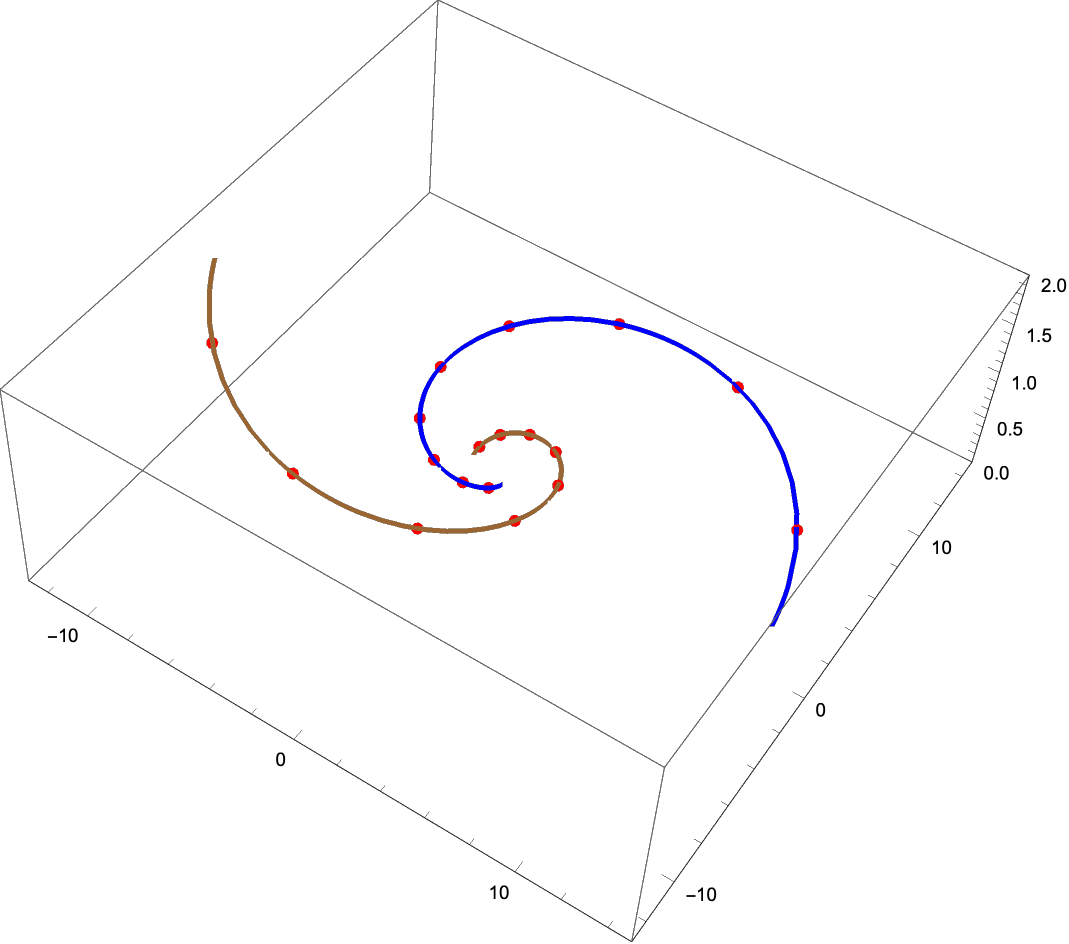} ~ \includegraphics[width=65mm]{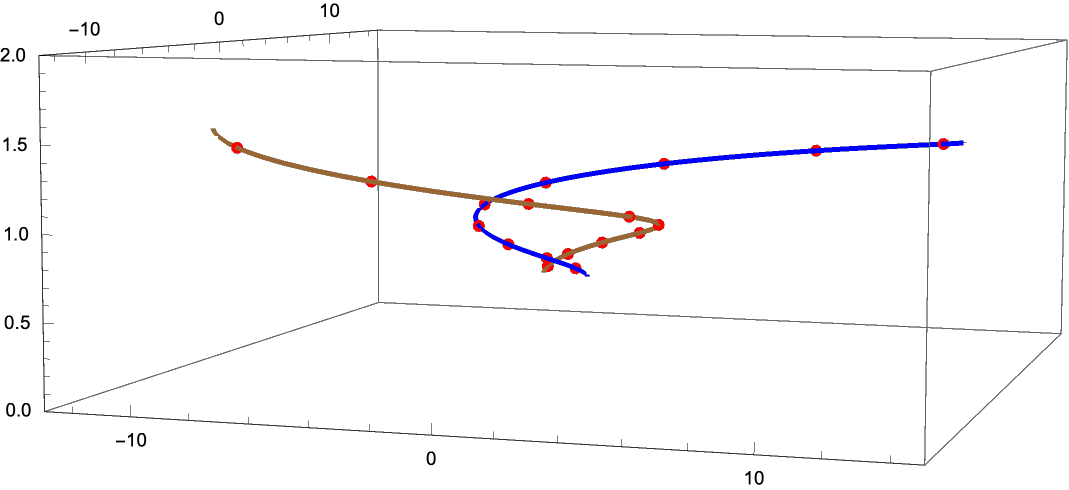}
\caption{Orbit (\ref{rotation3}) of  $(0.93, 0.5, -0.38)$ from different perspectives.} 

~

\includegraphics[width=65mm]{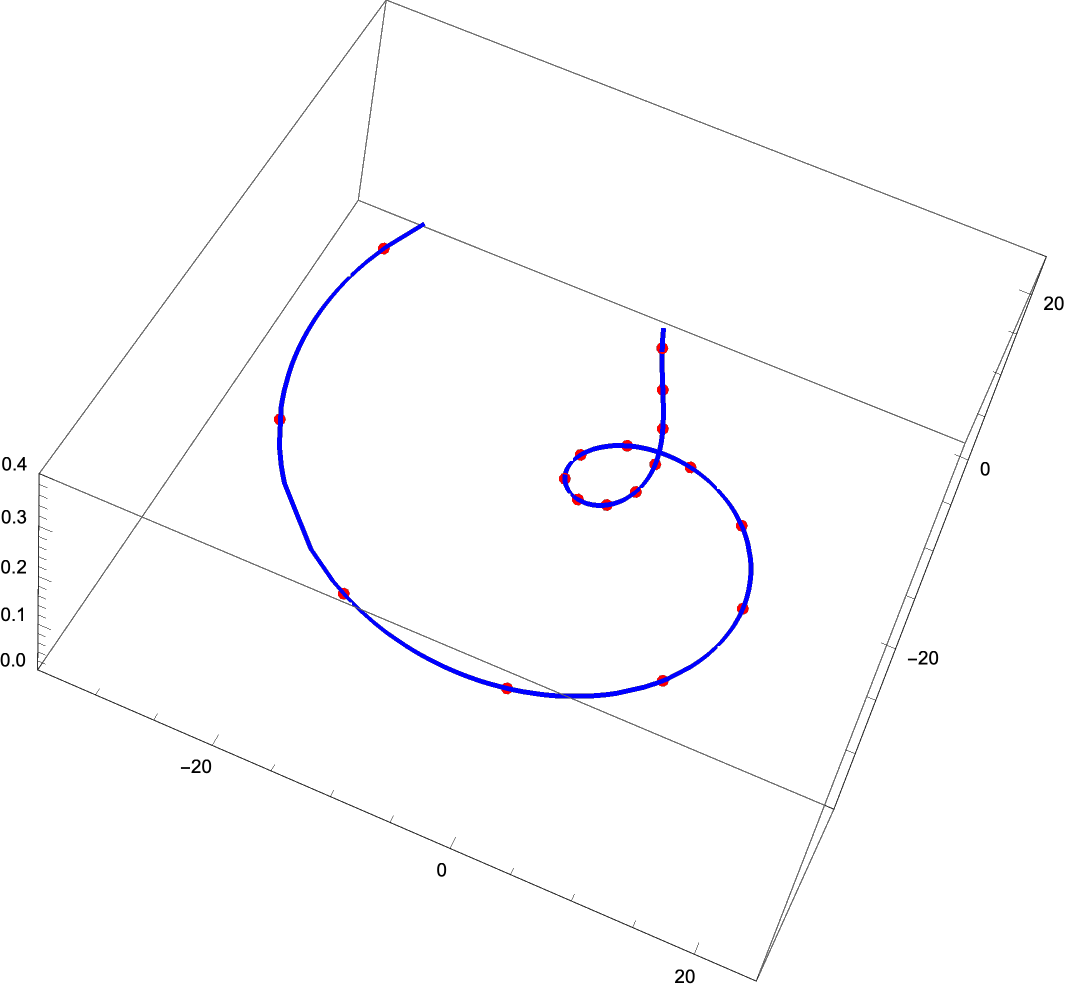} ~ \includegraphics[width=65mm]{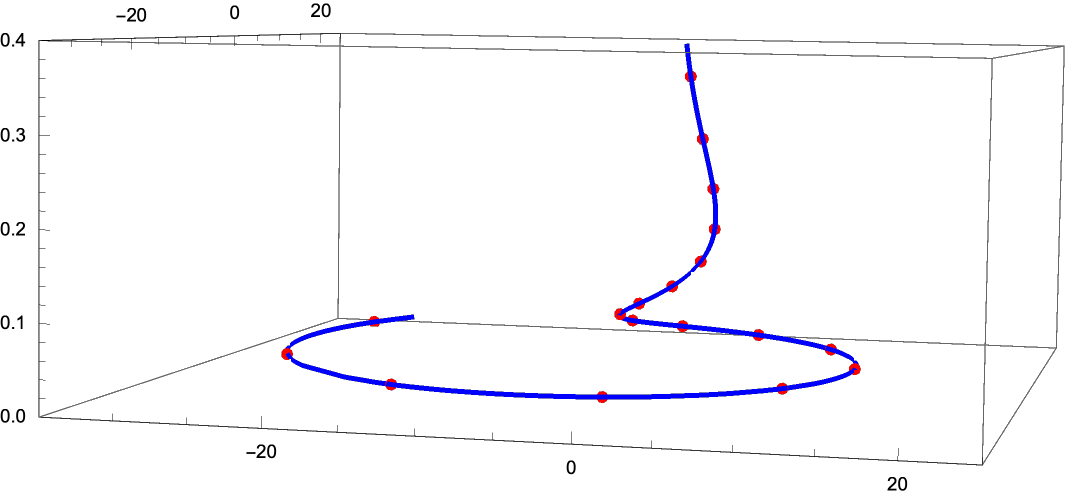}
\caption{Orbit (\ref{rotation3}) of  $(-0.5, 0.4, 0.89)$ from different perspectives.}

~

\centering
\includegraphics[width=65mm]{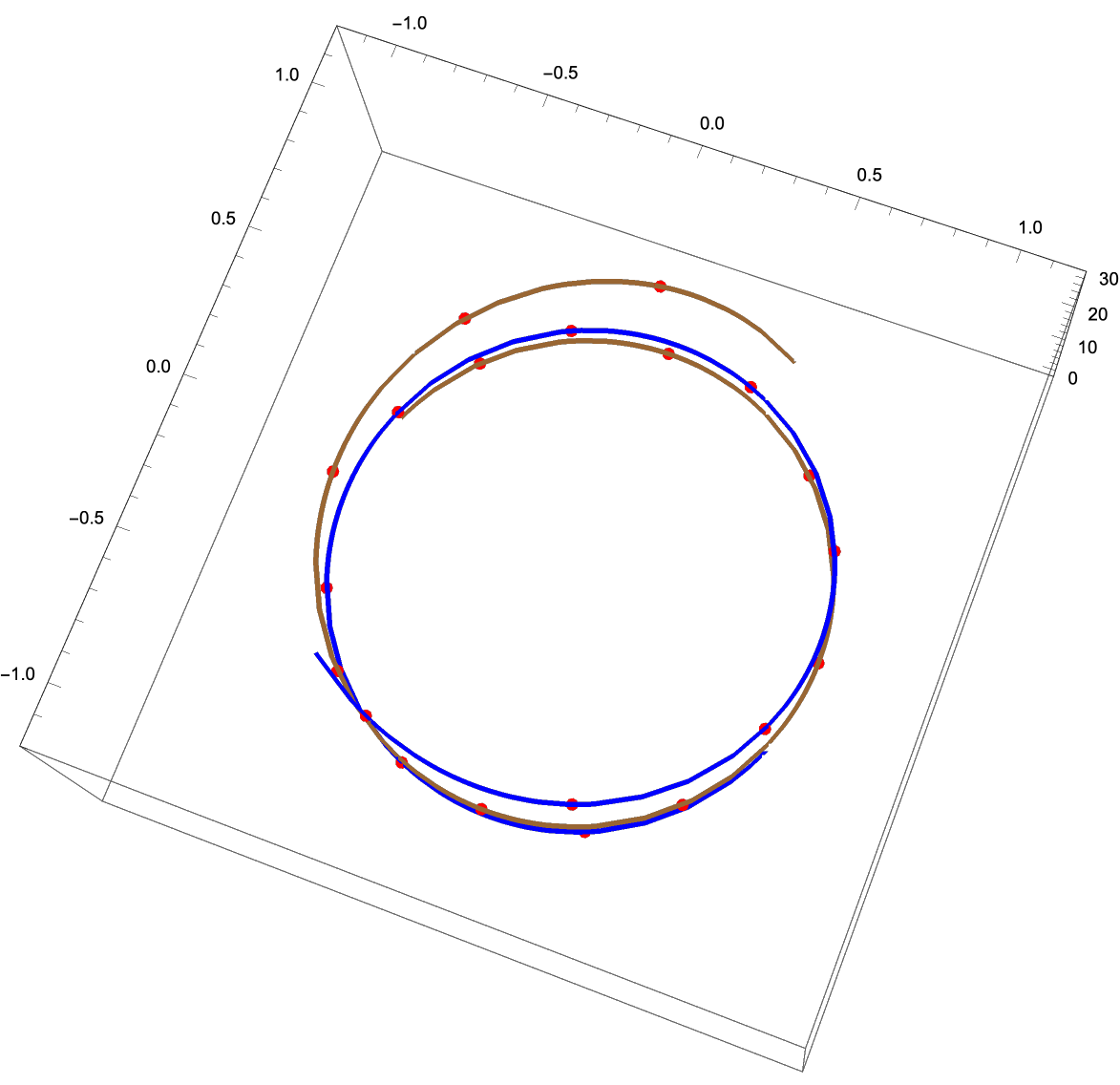} ~ \includegraphics[width=65mm]{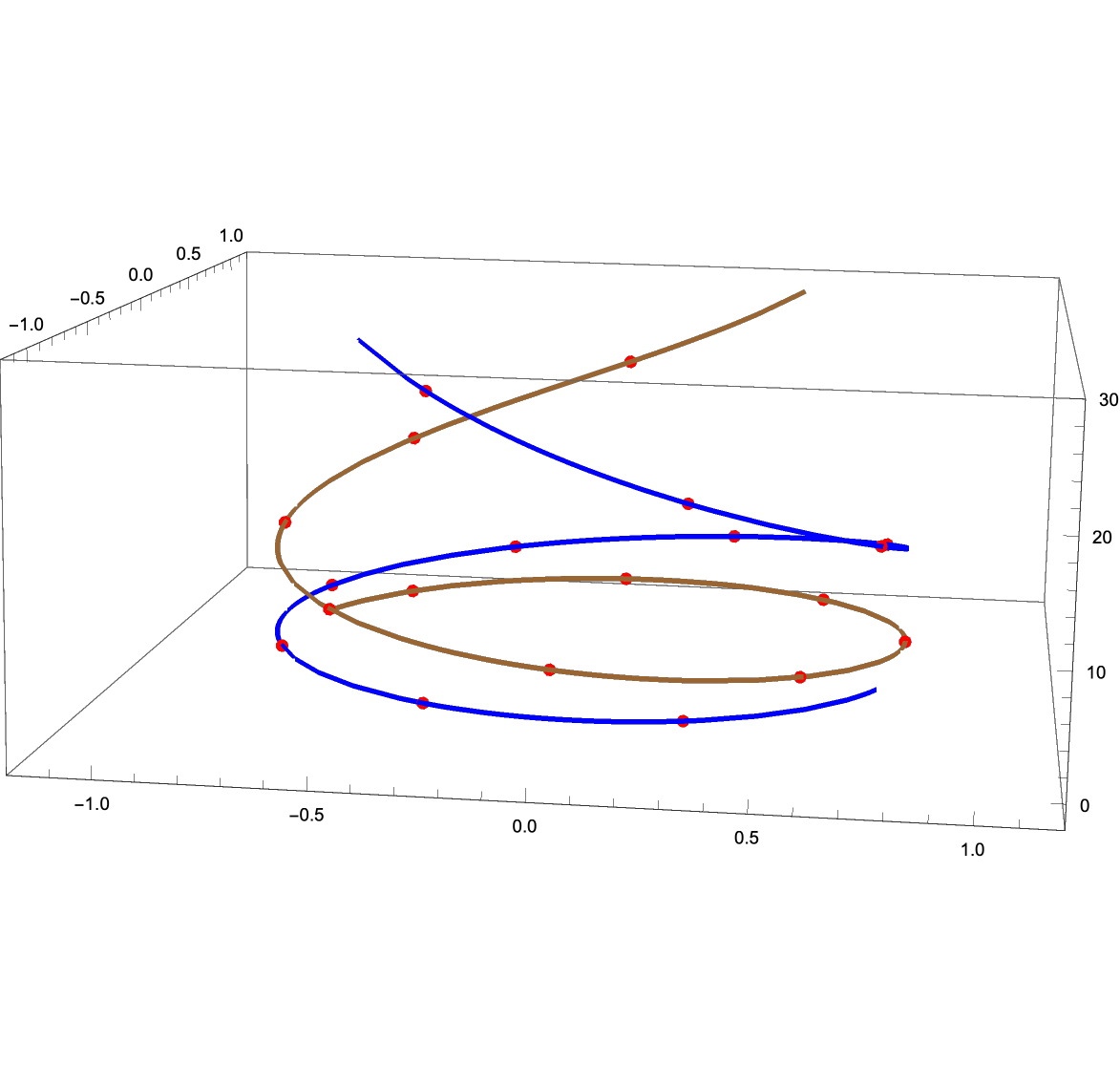}
\caption{Orbit (\ref{rotation3}) of  $(0.9289, 0.487, -0.2159)$ from different perspectives.}
\label{}
\end{figure}


\subsection{Examples of Horizontal Periodicity}

In this section I will discuss another type of periodicity, which the Riordan array of $p(t)$ 
may have. Such periodicity is naturally induced when $\widetilde{V}_{p(t)}$ has a finite order. 
It is clear from (\ref{rotation2}) that for linear polynomials to have a periodic orbit, we must have 
$(b - a)^{n+1} = b - a$ and $ (a + b)^{n+1} = a + b$, which over $\real$ either have the trivial 
solutions with $ab = 0$, two fixed points $\{(1/2,1/2),(-1/2,1/2)\}$, or two points with period 2:  
$\{(1/2,-1/2),(-1/2,-1/2)\}$. Consider for example, $p(t) = 1/2 - t/2$ with the following Riordan array:

\begin{equation}
\label{matrixVPer}
\left(\frac{1}{1-t^2}, \frac{t-t^2}{2}\right) = 
\left(
\begin{BMAT}[3pt]{ccccccc}{cccccccc}
1 & 0 & 0 & 0 & 0 & 0 &  0  \\
0 & 1/2 & 0  & 0 & 0 & 0 &  0  \\
1 & -1/2 & 1/4 &  0 & 0 & 0 &  0  \\
0 & 1/2 & -1/2 & 1/8 & 0 & 0 &  0  \\
1 & -1/2 & 1/2 & -3/8 & 1/16 & 0 &  0  \\
0 & 1/2 & -1/2 & 1/2 & -1/4 & 1/32 & 0 \\
1 & -1/2 & 1/2 & -1/2 & 7/16 & -5/32 &  1/64  \\
\vdots & \vdots &  \vdots & \vdots & \vdots & \vdots & \ddots 
\addpath{(1,1,1)ruuldd}
\addpath{(3,1,1)ruuldd}
\end{BMAT}
\right)
\end{equation}

We see two identical blocks in rows 5 and 6, and columns 1 and 3 as expected, but it is also interesting to look at the 
numbers in each column before the periodic cycles $(1/2, -1/2)$ and $(-1/2, 1/2)$ appear. 
If we add all such numbers in each column, we will get the following values.
$$
\begin{tabular}{|| m{3cm} || m{0.4cm} || m{0.4cm} || m{0.4cm} || m{0.4cm} || m{0.4cm} || m{0.4cm} ||  m{0.4cm} || m{0.4cm} ||}
\hline
\multicolumn{9}{||c||}{$p(t) = 1/2 - t/2$}\\ [3ex] 
\hline\hline
Column's \# k & 1 & 2 & 3 & 4 & 5 & 6 & 7 & $\ldots$ \\
\hline\
k-th column sum before the periodicity  &  0 & $\frac{1}{4}$ & $\frac{-1}{4}$ 
& $\frac{1}{4}$ & $\frac{-1}{4}$ & $\frac{1}{4}$ & $\frac{-1}{4}$ & $\ldots$ \\
\hline
\end{tabular}
$$
It seems that after the first column, such sums will repeat periodically with the values of 1/4 and -1/4. We will 
see below a similar behavior for the quadratic $p(t) = (-1 + 2t + 2t^2)/3$. I will prove the 
corresponding statement for the quadratic in Proposition 5 but leave it to the reader to show that the 
sum of the first $1 + 2(k-1)$ numbers in the $k$-th column of (\ref{matrixVPer}) 
is always $(-1)^k/4$. Notice that the table for $p(t) = -1/2 -t/2$ shows that there is no 
such periodicity of column sums in general.
$$
\begin{tabular}{|| m{3cm} || m{0.4cm} || m{0.4cm} || m{0.4cm} || m{0.4cm} || m{0.4cm} || m{0.4cm} ||  m{0.4cm} || m{0.4cm} ||}
\hline
\multicolumn{9}{||c||}{$p(t) = -1/2 - t/2$}\\ [3ex] 
\hline\hline
Column's \# k & 1 & 2 & 3 & 4 & 5 & 6 & 7 & $\ldots$ \\
\hline\
k-th column sum before the periodicity  &  0 & $\frac{1}{4}$ & $\frac{-1}{2}$ 
& $\frac{3}{4}$ & -1 & $\frac{5}{4}$ & $\frac{-3}{2}$ & $\ldots$ \\
\hline
\end{tabular}
$$
Let us now look at some quadratic polynomials with periodic orbits. With quadratics we can have many 
more such examples, and formula (\ref{rotation3}) implies that when  
$$
a + b + c = a^2 + b^2 + c^2 - ab - ac - bc = 1 ~~ \mbox{and} ~~ \cos\left(\frac{m\pi}{n}\right) = \frac{2 - 3(a + b)}{2}
$$
the orbit in the rotated coordinates will consist of finitely many points on a circle in the hyperplane 
$\tilde{z} = 1/\sqrt{3}$. Therefore, after we drop the first finitely many terms in each column of $RAp$, 
the corresponding infinite tails of the columns must repeat periodically vertically and horizontally. 
Consider for example the case when
$$
\cos\left(\frac{2\pi}{6}\right) = \frac{2 - 3(a + b)}{2}, ~ c = 1 - a - b,~ (a+b)^2 - (a+b) - ab = 0
$$ 
with solutions
$$
(a_1,b_1,c_1) = \bigl( -1/3, 2/3, 2/3\bigr) ~~ \mbox{and} ~~ (a_2,b_2,c_2) = \bigl(2/3, -1/3, 2/3\bigr).
$$
Take $p(t) = (-1/3) + (2/3)t + (2/3)t^2$ and write its $RAp$ in an ``abbreviated form" by forgetting all the terms 
that don't belong to the vertical periodic cycles. In other words, if we drop in the $k$th column (excluding $k=0$)  
first $3(k-1)$ numbers before the periodic cycle appears, 
and write the number of dropped terms in the $k$-th position of the zeroth (blue) row, we obtain the 
following doubly periodic infinite matrix 
\begin{equation}
\label{abbreviated}
``RAp" = \begin{bmatrix}
\rowcolor{blue!25}
0 & 1 & 4 & 7 & 10 & 13  & 16 & 19 & \cdots \\
1 & -1/3 & 0 & 2/3 & 1 & 2/3 & 0 & -1/3 & \cdots \\
0 & 2/3 & 1 & 2/3 & 0 & -1/3  & 0 & 2/3 & \cdots \\
0 & 2/3 & 0 & -1/3 & 0 & 2/3 & 1 & 2/3 & \cdots \\
1 & -1/3 & 0 & 2/3 & 1 & 2/3  & 0 & -1/3 & \cdots \\
0 & 2/3 & 1 &  2/3 & 0 & -1/3  & 0 & 2/3 & \cdots \\
0 & 2/3 & 0 &  -1/3 & 0 & 2/3 & 1 & 2/3 & \cdots \\
1 & -1/3 & 0 & 2/3 & 1 & 2/3  & 0 & -1/3 & \cdots \\
\vdots & \vdots & \vdots & \vdots & \vdots & \vdots & \vdots & \vdots & \ddots \\
\end{bmatrix}.
\end{equation}
Compare it with the actual Riordan array for $p(t) = (-1 + 2t + 2t^2)/3$.
\begin{equation}
\label{actual_RAp}
\begin{bmatrix}
1 & 0 & 0 & 0 & 0 & 0 &  \cdots \\
0 & -1/3 & 0  & 0 & 0 & 0 &  \cdots \\
0 & 2/3 & 1/9 &  0 & 0 & 0 &  \cdots \\
1 & 2/3 & -4/9 & -1/27 & 0 & 0 &  \cdots \\
0 & -1/3 & 0 & 2/9 & 1/81 & 0 &  \cdots \\
0 & 2/3 & 1 & -2/9 & -8/81 & -1/243 &  \cdots \\
1 & 2/3 & 0 & -17/27 & 16/81 & 10/243 &  \cdots \\
0 & -1/3 & 0 & 2/3 & 17/81 & -10/81 &  \cdots \\
0 & 2/3 & 1 & 2/3 & -64/81 & -1/243 & \cdots\\
1 & 2/3 & 0 & -1/3 & -16/81 & 130/243 & \cdots \\
\vdots & \vdots & \vdots & \vdots & \vdots & \vdots & \ddots \\
\end{bmatrix}
\end{equation}

Notice that even though the periodicity cycle $(1,0,0)$ of the zeroth column in (\ref{abbreviated}) 
coincides with the one in the fourth column, this is not the right point of view here. 
The zeroth column is an exception, and the right way to look at the periodicity {\sl horizontally
for columns}, is to start with the first column having the cycle $(-1/3, 2/3, 2/3)$. Exactly the same cycle 
appears the next time in the 7th column, and so on. 
As we see in this example, the {\sl horizontal prime period} for columns is 6, 
and the {\sl vertical prime period} for rows is 3. These numbers correspond to the rotation angle $2\pi/6$ 
and the degree 3 of $tp(t)$. 
 
Analogously to the linear polynomials $p(t) = (\pm 1)/2 - t/2$, one may ask, what would be the total 
sum of the dropped elements in each column for a quadratic polynomial? 
Checking first several columns in (\ref{actual_RAp}) we see again a periodicity with the period 6. 

$$
\begin{tabular}{|| m{3cm} || m{0.4cm} || m{0.4cm} || m{0.4cm} || m{0.4cm} || m{0.4cm} 
|| m{0.4cm} || m{0.4cm} || m{0.4cm} || m{0.4cm} || m{0.4cm} || m{0.3cm} ||}
\hline
\multicolumn{12}{||c||}{$p(t) = (-1/3) + (2/3)t + (2/3)t^2$}\\ [4ex] 
\hline\hline
Column's \# k & 1 & 2 & 3 & 4 & 5 & 6 & 7 & 8 & 9 & 10 & $\ldots$ \\
\hline\
k-th column sum before the periodicity  &  0 & $\frac{-1}{3}$ & $\frac{-2}{3}$ 
& $\frac{-2}{3}$ &$\frac{-1}{3}$ & 0 & 0 & $\frac{-1}{3}$  & $\frac{-2}{3}$ & $\frac{-2}{3}$ & $\ldots$ \\
\hline
\end{tabular} 
$$

\noindent However this is not true in general. Below is the corresponding table for the polynomial 
$p(t) = 2/3 - t/3 + 2t^2/3$, which also has $V_{p(t)}$ of order 6.

$$
\begin{tabular}{|| m{3cm} || m{0.4cm} || m{0.4cm} || m{0.4cm} || m{0.4cm} || m{0.4cm} 
|| m{0.4cm} || m{0.4cm} || m{0.4cm} || m{0.4cm} || m{0.4cm} || m{0.3cm} ||}
\hline
\multicolumn{12}{||c||}{$p(t) = (2/3) - (1/3)t + (2/3)t^2$}\\ [4ex] 
\hline\hline
Column's \# k & 1 & 2 & 3 & 4 & 5 & 6 & 7 & 8 & 9 & 10 & $\ldots$ \\
\hline\
k-th column sum before the periodicity  &  0 & 0 &$\frac{1}{3}$
& 1 &$\frac{5}{3}$ & 2 & 2 & 2  & $\frac{7}{3}$ & 3 & $\ldots$ \\
\hline
\end{tabular} 
$$

Now let us prove the periodicity of column sums for the polynomial $p(t) = -1/3 + 2t/3 + 2t^2/3$. 
The notations are the same as in \S 2.

\begin{prop} For $p(t) = (-1 + 2t + 2t^2)/3$, the circulant matrix $V_{p(t)}$ has order $6$, and for 
each $k\in \{1,2,\ldots,6\}$ and $\forall n\in\natu$,
\begin{equation}
\label{columns}
\sum\limits_{i = 0}^{3(k - 1)} C_{i,k} = \sum\limits_{i = 0}^{3(k + 6n - 1)} C_{i,k + 6n}.
\end{equation} 
\end{prop}
\begin{proof}
We leave it to the reader to check that $V^6_{p(t)} = I$. Recall that according to our Theorems 1 and 2, 
the periodicity in $k$-th column begins with the $1 + 3(k - 1)$-st term, and the periodic cycle in that column is given by 
$$
\begin{pmatrix}
a_{k-1} \\ b_{k-1} \\ c_{k-1} \\
\end{pmatrix} : = V^{k-1}_{p(t)} 
\begin{pmatrix}
-1/3 \\ 2/3 \\ 2/3 \\
\end{pmatrix}, ~ \mbox{where} ~ k \geq 1, ~ \mbox{and} ~ 
\begin{pmatrix}
a_0 \\ b_0 \\ c_0 \\
\end{pmatrix} \stackrel{def}{=}
\begin{pmatrix}
-1/3 \\ 2/3 \\ 2/3 \\
\end{pmatrix}.
$$
Hence we can represent the periodic part of the $k$th column by the generating function
$$ 
t^{1 + 3(k-1)} \cdot \frac{a_{k-1} + tb_{k-1} + t^2c_{k-1}}{1-t^3},
$$
and the sum of the first $3(k - 1)$ terms will be equal 
$$ 
\sum\limits_{i = 0}^{3(k - 1)} C_{i,k} t^i = \frac{(tp(t))^k}{1-t^3} - 
t^{1 + 3(k-1)}\cdot \frac{a_{k-1} + tb_{k-1} + t^2c_{k-1}}{1-t^3}.
$$
Similarly for the $(k+6n)$-th column 
$$
\sum\limits_{i = 0}^{3(k + 6n - 1)} C_{i,k+6n} t^i 
= \frac{(tp(t))^{k+6n}}{1-t^3} - 
t^{1 + 3(k+ 6n - 1)}\frac{a_{k + 6n-1} + tb_{k + 6n -1} + t^2c_{k + 6n-1}}{1-t^3}.
$$
Since $V_{p(t)}$ has order 6, we have 
$$
a_{k + 6n-1} + tb_{k + 6n -1} + t^2c_{k + 6n-1} = a_{k -1} + tb_{k -1} + t^2c_{k -1},
$$
and it is enough to show that
$$
\left[\frac{(tp(t))^k}{1-t^3} - 
t^{1 + 3(k-1)}\frac{a_{k-1} + tb_{k-1} + t^2c_{k-1}}{1-t^3} \right.
$$

\begin{equation}
\label{sumperiod}
-\left. \left(\frac{(tp(t))^{k+6n}}{1-t^3} - 
t^{1 + 3(k+ 6n - 1)}\frac{a_{k -1} + tb_{k -1} + t^2c_{k -1}}{1-t^3}\right)\right]_{t = 1} = 0.
\end{equation}

Let us abbreviate $a_{k -1} + tb_{k -1} + t^2c_{k -1}$ by $p_{k-1}(t)$ and simplify the difference  
$$
t^{1 + 3(k+ 6n - 1)}\frac{a_{k -1} + tb_{k -1} + t^2c_{k -1}}{1-t^3} - 
t^{1 + 3(k-1)}\frac{a_{k-1} + tb_{k-1} + t^2c_{k-1}}{1-t^3}
$$

$$
= t^{1 + 3(k-1)}\frac{p_{k-1}(t)}{1-t^3}\left(t^{3\cdot 6n} - 1\right) = 
t^{1 + 3(k-1)}\frac{p_{k-1}(t)}{1-t^3}(t^3 - 1)(1 + t^3 + \cdots + t^{3(6n-1)})
$$

\begin{equation}
\label{lastsum1}
= - t^{1 + 3(k-1)}(p_{k-1}(t))(1 + t^3 + \cdots + t^{3(6n-1)}).
\end{equation}
Our assumptions guarantee that $p_{k-1}(1) = 1, \forall k\geq 1$, and hence when $t = 1$, the product 
(\ref{lastsum1}) equals $-6n$. Thus, our goal written as (\ref{sumperiod}), now simplifies to the equality
\begin{equation}
\label{lastsum2}
\left[\frac{(tp(t))^k}{1-t^3} - \frac{(tp(t))^{k+6n}}{1-t^3} - 6n\right]_{t = 1} = 0.
\end{equation}
Since
$1 - tp(t) = 1 + (1/3)t - (2/3)t^2 - (2/3)t^3 = (1 - t)(3+4t+2t^2)/3$
we have
$$
\frac{(tp(t))^k - (tp(t))^{k+6n}}{1-t^3}  = \frac{(tp(t))^k\bigl(1- (tp(t))^{6n}\bigl)}{(1-t)(1+t+t^2)}
$$

$$
= (tp(t))^k\frac{(3+4t+2t^2)\bigl( 1 + tp(t) + \cdots (tp(t))^{6n-1}\bigr)}{3(1 + t + t^2)}.
$$
Hence, after substituting $t = 1$, the expression in the bracket in (\ref{lastsum2}) equals 
$$
\frac{9\overbrace{(1 + 1 + \cdots  + 1)}^{\mbox{$6n$ terms}}}{9} -  6n = 0
$$
as required.
\end{proof}

\noindent{Note:} It is the equality (\ref{lastsum2}), that doesn't hold true for a general case of a 
quadratic polynomial with a periodic orbit. The other steps of the proof can be easily generalized.


\section{Characterization by A- and Z- \\Sequences and related Combinatorics}

It is well known that each element of a  Riordan array can be described by a linear combination of the elements 
in the row above. The coefficients of such linear 
combinations do not depend on the row/column index and constitute two specific sequences. The Z-sequence 
is used for elements in the zeroth column, and the A-sequence for all other elements of the array. 
Initially such a description  was given by Rogers in \cite{Rogers}. In this last section we will assume that 
$R$ is a field of characteristic zero ($\real$ or $\comp$), and discuss the A- and Z- sequences for the 
Riordan array of $p(t) = a_0 + a_1t + \cdots + a_dt^d\in R[t]$. 
The reader will find all the details on formulas connecting these sequences to the Riordan array in sections 
5.2 and 5.3 of \cite{Barry} or chapter 4 of \cite{Shapiro2}. So, consider the array (\ref{array}) 
and, following notations of \cite{Shapiro2}, denote it by ${\R} (d(t),h(t))$. If we use $\bar{h}(t)$ for  
the compositional inverse of $h(t)$, i.e. $h(\bar{h}(t)) = \bar{h}(h(t)) = t$, then according to theorems 
4.3 and 4.5 of \cite{Shapiro2}, we get the generating functions of the Z and A-sequences correspondingly
\begin{equation}
\label{A-seq}
Z(t) = \frac{d(\bar{h}(t)) - d_{00}}{\bar{h}(t) d(\bar{h}(t))} ~~ \mbox{and} ~~ A(t) = \frac{t}{\bar{h}(t)},
\end{equation}
where $d_{00}$ is the element of ${\R} (d(t),h(t))$ in the 0-th row and 0-th column. Since in our cases  
$d(t) = 1/(1- t^{d+1}),~h(t) = tp(t)$ and $d_{00} = 1$, we can rewrite (\ref{A-seq}) as
\begin{equation}
\label{A-seq2}
Z(t) = \bar{h}^d(t) ~~ \mbox{and} ~~ A(t) = \frac{t}{\bar{h}(t)}.
\end{equation}
Notice that the f.p.s. $\bar{h}(t) = \alpha_1 t+\alpha_2 t^2 + \cdots$ is the solution of the equation
$$
\bar{h}(t)p(\bar{h}(t)) = t \Longleftrightarrow a_0\bar{h}(t) + a_1\bar{h}^2(t) + \cdots + a_d\bar{h}^{d+1}(t) = t.
$$
Now, if $p(t) = a$ then $\bar{h}(t) = t/a$ and 
$$
Z(t) = 1 ~~ \mbox{and} ~~ A(t) = a.
$$
If $p(t) = a + bt,~b\neq 0$ then $\bar{h}(t)$ is the solution of $a\bar{h}(t) + b\bar{h}^2(t) = t$, 
and since $\bar{h}(t)$ must have order 1, we have $\bar{h}(0) = 0$, so 
\begin{equation}
\label{Catalan1}
Z(t) = \bar{h}(t) = \frac{\sqrt{a^2+4bt} - a}{2b} ~~ \mbox{and} ~~ A(t) = \frac{2bt}{\sqrt{a^2+4bt} - a}.
\end{equation}
In particular,
$$
Z(t) =  \left(\frac{-a}{b}\right)\frac{1 - \sqrt{1 - 4(-b/a^2)t} }{2} = 0t^0 + \left(\frac{-a}{b}\right) 
\sum\limits_{n = 0}^{\infty} C_n \left(\frac{-bt}{a^2}\right)^{n+1},
$$
where $C_n = \binom{2n}{n}/(n+1),~n\geq 0$ is the $n$th Catalan number (see for example 
\S 3.2 of \cite{Barry} or \S 1.2 of \cite{Shapiro2}). As for the A-sequence, since
$$
\frac{2bt}{\sqrt{a^2+4bt} - a}   = \frac{\sqrt{a^2+4bt} + a}{2},
$$
we similarly obtain 
$$
A(t) = \frac{\sqrt{a^2+4bt} + a}{2} = a\frac{1 + \sqrt{1 + 4(b/a^2)t} }{2} = a -  
a\sum\limits_{n = 0}^{\infty} C_n \left(\frac{-bt}{a^2}\right)^{n+1}.
$$
Consider next $p(t) = a + bt + ct^2$, so $\bar{h}(t)$ is the solution of a cubic equation 
$a\bar{h}(t) + b\bar{h}^2(t)  + c \bar{h}^3(t) = t$. The general formula in this case is already too 
complicated to shed some light, but we can easily obtain the first several terms 
of the f.p.s. expansion of $\bar{h}(t)$. Indeed, writing 
$$
\bar{h}(t) = \alpha_1t + \alpha_2 t^2 + \alpha_3 t^3 + \alpha_4 t^4 +\cdots, 
$$
substituting this sum into $\bar{h}(t) p(\bar{h}(t)) = t$, and comparing the coefficients of powers of $t$, 
one derives the following infinite system of equations.
$$
\left\{
\begin{array}{lcl}
a\alpha_1 = 1\\
a\alpha_2 + b\alpha_1^2 = 0\\
a\alpha_3 + 2b\alpha_1\alpha_2 + c\alpha_1^3 = 0\\
a\alpha_4 + b(\alpha_2^2 + 2\alpha_1\alpha_3) + 3c\alpha_1^2\alpha_2 = 0\\
~~~~~~~~~~~~~~~~~~~ \cdots\\
\end{array}
\right.
$$
Solving this system recursively we obtain
$$
\alpha_1 = \frac{1}{a},~ \alpha_2 = -\frac{b}{a^3}, ~ \alpha_3 = \frac{2b^2 - ac}{a^5}, ~ 
\alpha_4 = \frac{5 (-b^3 + a b c)}{a^7}, ~ \ldots.
$$
Therefore
$$
Z(t) = \bar{h}^2(t) = \frac{1}{a^2}t^2 - \frac{2b}{a^4}t^3 + \frac{5b^2 - 2ac}{a^6}t^4 + \cdots
$$
and 
$$
A(t) = a + \frac{b}{a}t + \frac{-b^2 + ac}{a^3}t^2 + \frac{2b^3-3abc}{a^5}t^3 + 
\frac{-5b^4 + 10a b^2 c - 2a^2 c^2}{a^7}t^4 + \cdots.
$$
Notice that the sequences of the numerical coefficients of the powers of $b$ in the last 
expansions for $Z(t)$ and $A(t)$ contain Catalan numbers (up to a sign) as well:
$$
1,-2,5,\ldots ~~ \mbox{and} ~~ 0,1,-1,2,-5, \ldots.
$$
It is not a surprise, since if we let $c=0$ in the Riordan array $RAp$, we will get the array with 
$h(t) = t(a+bt)$, which has the compositional inverse $\bar{h}(t)$ given in (\ref{Catalan1}). The same, of course,  
will be true for a polynomial $p(t)$ of any degree $d\geq 1$. 

It is interesting to note that not only the coefficients of the powers of $b$, but also the 
numerical coefficients of the terms containing $ct^n,~ n\geq 0$ in the expansion of $A(t)$ 
have some combinatorial interpretation. To simplify the formula of $A(t)$, 
take for instance $p(t) = 1 + t + ct^2$. Then with some help of a computer algebra system 
(e.g. {\sl Mathematica}) one finds that for this particular $p(t)$ we have
$$
\bar{h}(t) = t - t^2 + (2 - c) t^3 + 5 (-1 + c) t^4 + (14 - 21 c + 3 c^2) t^5 + \cdots
$$
and the formal power series $A(t)$ begins with 
$$
1 + t + (-1+c)t^2 + (2-3c)t^3 + (-5 + 10 c - 2 c^2)t^4 + (14 - 35 c + 15 c^2)t^5 
$$
\begin{equation}
\label{c_sum}
 + (-42 + 126 c  - 84 c^2 + 7 c^3)t^6 + (132 - 462 c + 420 c^2 - 84 c^3) t^7 + \cdots 
\end{equation}
The number coefficients of $(c)t^n$ in (\ref{c_sum}) are $0,0,1,-3,10,-35,126,-462$. We can also 
obtain these numbers as binomial coefficients (a proof is given below)
\begin{equation}
\label{binom2}
\left\{(-1)^n\binom{2n+1}{n+1} \right\}_{n=0}^{\infty} = \{1, -3, 10, -35, 126, -462, \ldots\}.
\end{equation}
Absolute values of these numbers give the total number of leaves in all rooted ordered trees with $n\geq 0$ 
edges (see sequences \href{https://oeis.org/A001700/}{A001700} and 
\href{https://oeis.org/A088218/}{A088218} in the \cite{OEIS}, the On-Line Encyclopedia of Integer Sequences).

Moreover, the numerical coefficients of the terms containing $(c^2)t^n,~ n\geq 0$ in (\ref{c_sum}) are 
$\{0,0,0,0,-2,15,-84,420,\ldots\}$. For a combinatorial description of these numbers 
in terms of Dyck words, we refer the reader to 
\href{https://oeis.org/A002740/}{A002740}. Now let me show  
that the numerical coefficients of $(a^ib^jc)t^n,~ n\geq 0$ in $A(t)$ for any quadratic 
$p(t) = a+bt+ct^2$ with $ac\neq 0$ are indeed given by (\ref{binom2}). 
To isolate such coefficients in the f.p.s. expansion of $A(t)$, it is enough to differentiate 
$A(t)$ once with respect to the parameter $c$, and then substitute $c = 0$. Hence our goal is
\begin{theorem}
For any $n\in\natu$,
$$
[t^{n+2}]\left.\frac{d(A(t))}{dc}\right|_{c=0} = \frac{(-b)^n}{a^{2n+2}}\binom{2n+1}{n+1} 
$$
\end{theorem}
\begin{proof}
Clearly we can first apply the {\sl ``coefficient of operator"}, and then differentiation with substitution. Since 
$A(t) = t/\bar{h}(t)$, using shifting together with the Lagrange inversion formula we can write $[t^n]A(t)$ as 
$$
[t^n]\frac{t}{\bar{h}(t)} = [t^{n-1}](\bar{h}(t))^{-1} = \frac{-1}{n-1}[t^n]\left(\frac{t}{h(t)}\right)^{n-1} = 
 \frac{-1}{n-1}[t^n]\left(\frac{1}{p(t)}\right)^{n-1}.
$$
Since $p(t) = a+ bt+ct^2$, let me write $p(t) = L + Qc$ where $L = a+ bt$ and $Q=t^2$. Then 
the powers of $1/p(t)$ will equal  
\begin{equation}
\label{sum4}
\left(\frac{1}{p(t)}\right)^{n-1} = \left(\frac{1}{L + Qc}\right)^{n-1} = 
\frac{1}{L^{n-1}}\left(1 - \frac{Q}{L}c + \left(\frac{Q}{L}c\right)^2 + \cdots \right)^{n-1}.
\end{equation}
Since all the terms involving $c^k,~k\geq 2$ will turn to zero after differentiating and substituting 
$c=0$, we can deduce from (\ref{sum4}) that 
$$
[t^n]\left.\frac{d(A(t))}{dc}\right|_{c=0} = \frac{-1}{n-1}[t^n]\left.\frac{d}{dc}\right |_{c=0} 
\frac{1}{L^{n-1}}\left(1 - \frac{Q}{L}c\right)^{n-1}
$$
\begin{equation}
\label{last}
=  \frac{-1}{n-1}[t^n]\left.\frac{n-1}{L^{n-1}}\left(1 - \frac{Q}{L}c\right)^{n-2}\frac{-Q}{L}\right|_{c=0} 
= [t^n]\frac{Q}{L^n} = [t^n]\frac{t^2}{(a + bt)^n}.
\end{equation}
Simplifying and applying the so-called Newton's rule (see \S 2.2. of \cite{Shapiro2}), 
we continue (\ref{last}) with 
$$
[t^n]\frac{t^2}{(a + bt)^n} = \frac{1}{a^n} [t^{n-2}]\bigl(1 + (b/a)t\bigr)^{-n} = 
a^{-n}\binom{-n}{n-2}\left(\frac{b}{a}\right)^{n-2}
$$

$$
= \frac{(-b)^{n-2}}{a^{2n-2}}\binom{n + n-2 -1 }{n-2} = 
\frac{(-b)^{n-2}}{a^{2n-2}}\binom{2n -3 }{n-2}.
$$
Finally, replacing $n$ by $n+2$, we see that
$$
[t^{n+2}]\left.\frac{d(A(t))}{dc}\right|_{c=0} = \frac{(-b)^n}{a^{2n+2}}\binom{2n + 1}{n} 
= \frac{(-b)^n}{a^{2n+2}}\binom{2n + 1}{n+1}.
$$
\end{proof}

\end{document}